\documentclass[a4paper, 12pt]{article}
\usepackage{amssymb}
\usepackage{amsmath}
\usepackage{amsthm}
\usepackage[dvips]{graphicx}
\usepackage{verbatim,enumerate}
\usepackage{amssymb}
\theoremstyle{definition}
 \newtheorem{theorem}{Theorem}[section]
 \newtheorem*{theorem*}{Theorem}
 \newtheorem*{lemma*}{Lemma}
 
 \newtheorem{fact}[theorem]{Fact}
 \newtheorem{fact*}{Fact}
 \newtheorem{lemma}[theorem]{Lemma}

 \newtheorem*{remark*}{Remark}
 
 \newtheorem{example}[theorem]{Example}

\newcommand{\vect}[1]{\boldsymbol{#1}}
\newcommand{\R}{\boldsymbol{R}}
\newcommand{\trace}{\operatorname{trace}}
\newcommand{\rank}{\operatorname{rank}}
\newcommand{\Hess}{\operatorname{Hess}}
\renewcommand{\phi}{\varphi}

\newcommand{\ep}{\varepsilon}
\newcommand{\x}{\vect{x}}
\renewcommand{\u}{\vect{u}}
\renewcommand{\dfrac}{\displaystyle\frac}
\newcommand{\pmt}[1]{{\begin{pmatrix} #1  \end{pmatrix}}}
\begin{document}
\begin{center}
{\large 
{\bf Criteria for singularities of smooth maps from the plane into the
 plane and their applications}}
\\ Kentaro Saji
\\ \today
\end{center}
\renewcommand{\abstractname}{}
\begin{abstract}
 We give useful criteria of lips, beaks and swallowtail
 singularities of a smooth map from the plane into the plane.
 As an application of criteria, we discuss
 the singularities of a Cauchy problem of
 single conservation law.
\end{abstract}
\renewcommand{\thefootnote}{\fnsymbol{footnote}}
\footnote[0]{2000 Mathematics subject classification. 57R45}
 \section{Introduction}
\begin{figure}[ht]
\centering
\includegraphics[width=0.2\linewidth]{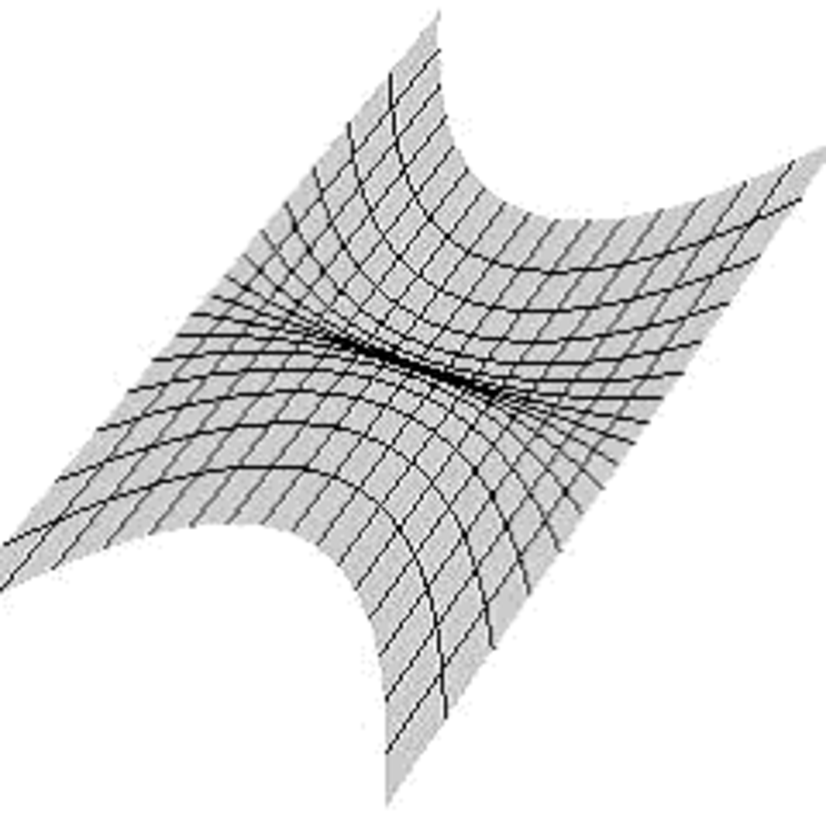}\hspace{10mm}
\includegraphics[width=0.2\linewidth]{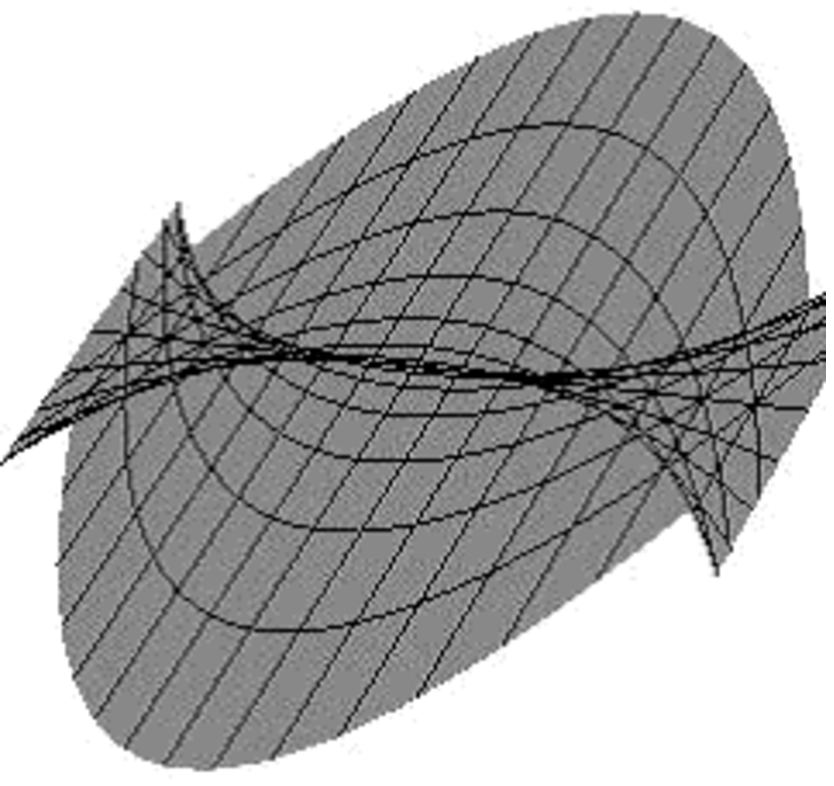}\hspace{10mm}
\includegraphics[width=0.2\linewidth]{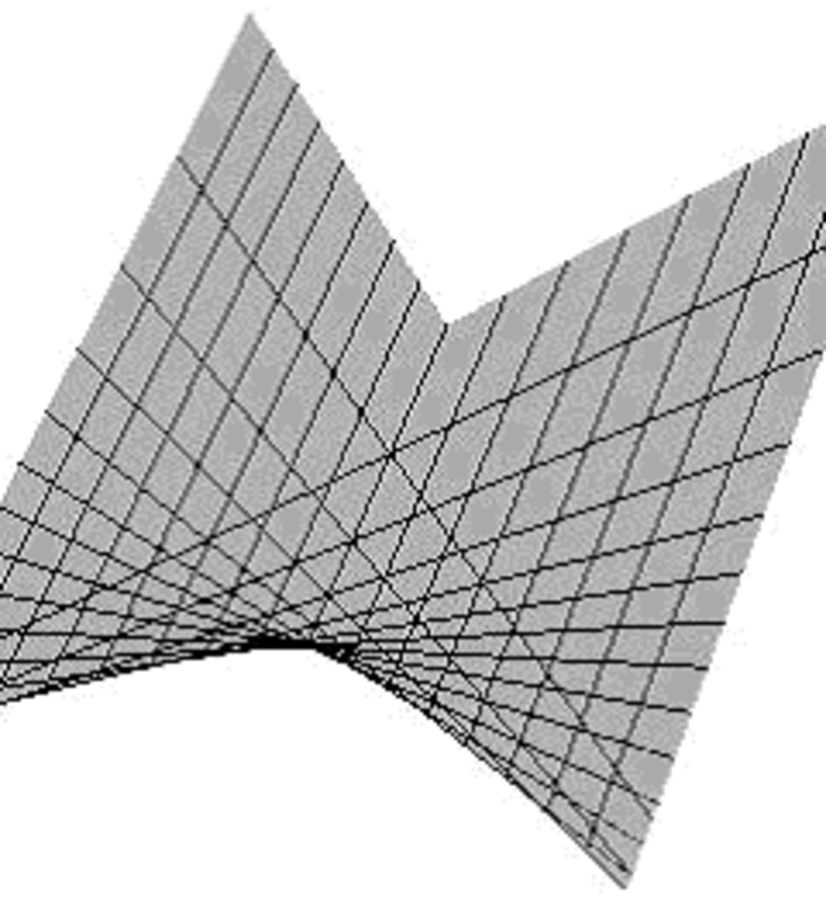}
\caption{Lips, beaks and swallowtail}
\label{fig:lipetc}
\end{figure}

Singularities of map germs have long been studied,
especially up to the equivalence
under coordinate changes
in both source and target (${\mathcal A}$-equivalence).
According to \cite{gaff}, ``classification'' 
for map germs with ${\mathcal A}$-equivalence
means finding lists of germs, and showing that
all germs satisfying certain conditions are equivalent to a germ
on the list.
Classification is well understood, with many good references in
the literature.
``Recognition'' means 
finding criteria which will describe
which germ on the list a given germ is equivalent to (see \cite{gaff}).
The classification problem and recognition problem for
map germs from the plane into the plane up to
${\mathcal A}$-equivalence
was studied by J. H. Rieger \cite{rie}.
He classified map germs $(\R^2,0)\to(\R^2,0)$
with corank one and ${\mathcal A}_e$-codimension $\leq 6$.
Table \ref{tab:rie}
shows the list of the ${\mathcal A}_e$-codimension $\leq3$
local singularities obtained in \cite{rie}.
 Some of these singularities are
also called as follows: 
$4_{2,+}$ ({\it lips\/}), $4_{2,-}$ ({\it beaks\/}),
$5$ ({\it swallowtail\/}). 
These singularities are depicted in Figure \ref{fig:lipetc}.
Rieger also discussed the recognition of these map germs
after normalizing the coordinate system as
$(u,v)\mapsto(u,f_2(u,v))$.
However, for applications, 
criteria of recognition without using
normalization are not only more convenient but also indispensable in
some cases.
\begin{table}[htbp]
 \begin{center}
  \begin{tabular}{llc}
   \hline
   Name&Normal form&${\mathcal A}_e$-codimension\\
   \hline
   Immersion&$(u,v)$&0\\
   Fold     &$(u,v^2)$&0\\
   Cusp     &$(u,v^3+uv)$&0\\
   $4_{k,\pm}$    &$(u,v^3\pm u^kv),\ k=2,3$&$k-1$\\
   $5$      &$(u,uv+v^4)$&1\\
   $6_\pm$      &$(u,uv+v^5\pm v^7)$&2\\
   $11_5$   &$(u,uv^2+v^4+v^5)$&2\\
   \hline
  \end{tabular}
  \label{tab:rie}
  \caption{Classification of $(\R^2,0)\to(\R^2,0)$}
 \end{center}
\end{table}

In this paper, we give criteria
for the lips, the beaks and the swallowtails of
a map germ $(\R^2,0)\to(\R^2,0)$
without using the normalizations (Theorem \ref{thm:main}).
Since they only use the information of the 
Taylor coefficients of the germ,
Theorem \ref{thm:main} 
can be applied directly for the recognition of
the lips, the beaks and the swallowtail
on explicitly parameterized maps.
Using the criteria, we study singularities of a
conservation law about a time variable.
We study singularities of geometric
solutions of the equation and show the singularities that
appear for the first time are generically the lips (Section 3).

The case of wave front surfaces in 3-space,
criteria for the cuspidal edge and the swallowtail
were given by M. Kokubu et al.\ \cite{krsuy}.
By using them, we studied
local and global behaviors of flat fronts in
hyperbolic 3-space.
Using them, 
K. Saji et al.\  \cite{front}
introduced the singular curvature
on the cuspidal edge and investigated its properties.

Criteria for other singularities of fronts 
and their applications
were given in \cite{fsuy,horo,suy3}.
Recently, several applications of
these criteria were considered in various situations
\cite{ishi-machi,horo,circular,kruy,sch}.

Throughout this paper,
we work in the $C^\infty$-category.

\section{Preliminaries and statements of criteria}
Let $U\subset\R^2$ be an open set and
$f:(U,p)\to(\R^2,0)$ a map germ.
We call $q\in U$ a {\it singular point\/} of $f$ if
$\rank(df)_q\leq1$.
We denote by $S(f)\subset U$ the set of singular points of $f$.
Two map germs $f_i:(\R^2,0)\to(\R^2,0)$ $(i=1,2)$ are
{\it ${\mathcal A}$-equivalent\/} if
there exist diffeomorphism map germs
$\Phi_i:(\R^2,0)\to(\R^2,0)$ $(i=1,2)$ such that
$f_1\circ\Phi_1=\Phi_2\circ f_2$ holds.
For a positive integer $k$,
 a map germ $f:(U,p)\to(\R^2,0)$ is $k$-{\it determined\/} if
any $g:(U,p)\to(\R^2,0)$ satisfying the condition that
the $k$-jet $j^kg(p)$ of $g$ is equal to $j^kf(p)$,
is ${\mathcal A}$-equivalent to $f$.
The following fact is well-known.
\begin{fact}$($\cite[Lemma 3.2.2 and 3.1.3]{rie}$)$
 \label{fact:det}
 The lips and the 
 beaks\/ $(x,y)\mapsto(x,y^3\pm xy)$ are three-determined.
 The swallowtail\/ $(x,y)\mapsto(x,xy+y^4)$ is four-determined.
\end{fact}

Let $f:(U,p)\to(\R^2,0)$ be a map germ.
A singular point $q$ is of {\it corank one\/} if
$\rank (df)_q=1$.
If $p$ is a corank one singular point of $f$,
then there exists a neighborhood $V$ of $p$ and a
never vanishing
vector field $\eta\in {\mathfrak X}(V)$ such that
$df_q(\eta)=0$ holds for any $q\in S(f)\cap V$.
We call $\eta$ the {\it null vector field}.
We define a function which plays a crucial role
in our criteria.
Let $(u_1,u_2)$ be coordinates of $U$.
Define the {\it discriminant function\/} $\lambda$ {\it of} $f$ by
\begin{equation}
 \label{eq:lambda}\lambda(u_1,u_2)=
  \det\left(\frac{\partial f}{\partial u_1},
       \frac{\partial f}{\partial u_2}\right)(u_1,u_2).
\end{equation}
Then $S(f)=\lambda^{-1}(0)$ holds.
We call $p\in S(f)$ a {\it non-degenerate singular point\/}
if
$d\lambda(p)\ne0$ and a {\it degenerate singular point\/}
if $d\lambda(p)=0$.
Note that a non-degenerate singular point is of corank one.
The terminologies ``discriminant function'',
``null vector field'' and ``non-degeneracy'' are
defined in \cite{krsuy}
in order to state criteria for fronts in the $3$-space.
Our definitions of these three terminologies are similar.
These notions also play a key role to identify
singularities for our case.
This seems to be related to the correspondence
between singularities of front and
its projection to the limiting tangent plane.
This correspondence is discussed in \cite{suy3}.
 
We review the criteria for the fold and the cusp,
due to Whitney \cite{whit} (see also \cite{suy3}).
\begin{fact}$($\cite[Proposition 2.1]{whit}$)$
\label{fact:whit}
 For a map germ\/ $f:(U,p)\to(\R^2,0)$,
 $f$ at\/ $p$ is\/ ${\mathcal A}$-equivalent to
 the fold if and only if\/
 $\eta\lambda(p)\ne0$.
 
 Furthermore, $f$ at\/ $p$ is\/ ${\mathcal A}$-equivalent to
 the cusp if and only if\/
 $p$ is non-degenerate, $\eta\lambda(p)=0$
 and\/ $\eta\eta\lambda(p)\ne0$.
\end{fact}
Here, $\eta\lambda$ means the directional derivative $D_\eta\lambda$.
The main result of this paper is the following.
\begin{theorem}
 \label{thm:main}
 For a map germ\/ $f:(U,p)\to(\R^2,0)$, the following hold.
 \begin{itemize}
  \item[(1)]  
              $f$ is\/ ${\mathcal A}$-equivalent to
              the lips if and only if\/ $p$ is of corank one,
              $d\lambda(p)=0$
              and\/ $\lambda$ has a Morse type
              critical point of index\/ $0$ or\/ $2$ at\/ $p$,
              namely, $\det\Hess\lambda(p)>0$.

  \item[(2)]$f$ is\/ ${\mathcal A}$-equivalent to
            the beaks if and only if\/ $p$ is of corank one,
            $d\lambda(p)=0$, $\lambda$ has a Morse type
            critical point of index\/ $1$ at\/ $p$ $($i.e.,
            $\det\Hess\lambda(p)<0$.$)$ and\/
            $\eta\eta\lambda(p)\ne0$.
            
  \item[(3)]$f$ is\/ ${\mathcal A}$-equivalent to
            the swallowtail if and only if\/
            $d\lambda(p)\ne0$,
            $\eta\lambda(p)=\eta\eta\lambda(p)=0$ and\/
            $\eta\eta\eta\lambda(p)\ne0$.
 \end{itemize}
\end{theorem}
Here, for a function $\lambda:(U,u_1,u_2)\to\R$,
$\Hess \lambda$ is the matrix defined by
$\Hess\lambda=(\partial^2\lambda/\partial u_i\,\partial u_j)_{i,j=1,2}$.
Remark that in Theorem \ref{thm:main} (1), $\eta\eta\lambda(p)\ne0$
is automatically satisfied because of the
symmetricity of $\Hess\lambda$ and
the inequality $\det\Hess\lambda(p)>0$.
\begin{example}
 Let us put
 $$
 \begin{array}{l}
 f_{\text{l}}(u,v)=(u,v^3+u^2v),\quad
 f_{\text{b}}(u,v)=(u,v^3-u^2v)\\
 \hspace{50mm}\text{and}\quad
 f_{\text{s}}(u,v)=(u,v^4+uv).
 \end{array}
 $$
 Since these are nothing but the defining formula for the lips, the beaks and
 the swallowtail, these maps satisfy the conditions in Theorem\/
 $\ref{thm:main}$.
 The discriminant functions for these maps are
 $$\lambda_{\text{l}}=3v^2+u^2,\quad
 \lambda_{\text{b}}=3v^2-u^2\quad\text{and}\quad
 \lambda_{\text{s}}=4v^3+u,
 $$
 respectively.
 Thus\/ $\lambda_{\text{l}}$ and\/ $\lambda_{\text{b}}$ have
 a Morse type critical point at the
 origin.
 Furthermore, the null vector field can be chosen as\/
 $\eta=(0,1)$ for all maps.
 It holds that\/ $\eta\eta\lambda_{\text{b}}\ne0$, and that\/
 $d\lambda_{\text{s}}\ne0$, $\eta\lambda_{\text{s}}=\eta\eta\lambda_{\text{s}}=0$
 and\/ $\eta\eta\eta\lambda_{\text{s}}\ne0$ at the origin.
 Thus we see that each of the conditions in Theorem\/ $\ref{thm:main}$ 
 is satisfied for each map.
 These observations together with the following 
 Lemma\/ $\ref{lem:indeplamb}$
 confirm the only if part of Theorem\/ $\ref{thm:main}$.
\end{example}
\begin{example}
 Let\/ $\gamma:I\to\R^2$ be a plane curve with\/ $\gamma'(t)\ne0$ 
 for any $t\in I$.
 The {\it tangential ruling map\/ $R_\gamma$ of\/ $\gamma$} is
 the map\/ $R_\gamma:(t,u)\mapsto \gamma(t)+u\gamma'(t)$.
 The discriminant function and the null vector field
 of\/ $R_\gamma$ are\/ $\lambda=u\kappa$ and\/ $\eta=(-1,1)$, respectively,
 where\/ $\kappa$ is the curvature of\/ $\gamma$.
 Thus we have
 $$\Hess\lambda(t,0)=\pmt{0&\kappa'\\ \kappa'&0}\quad\text{and}\quad
 \eta\eta\lambda(t,0)=-2\kappa'.
 $$
 Using Theorem\/ $\ref{thm:main}$,
 $R_\gamma$ at\/ $(t_0,0)$ is\/ ${\mathcal A}$-equivalent to
 the beaks if and only if\/ $\kappa(t_0)=0$ and\/ $\kappa'(t_0)\ne0$ holds\/
  $($See figure\/
  $\ref{fig:rulingex})$.
 \begin{figure}[ht]
\centering
\includegraphics[width=0.25\linewidth]{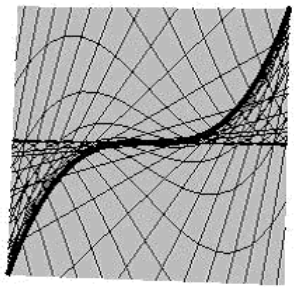}\\
\caption{The beaks on the tangential ruling map of $(t,t^3)$ at $(t,u)=(0,0)$.}
\label{fig:rulingex}
\end{figure}
\end{example}
To prove Theorem \ref{thm:main}, we need the following lemma.
\begin{lemma}
 \label{lem:indeplamb}
 For a map germ $f:(U,p)\to(\R^2,0)$,
 the conditions in Theorem\/ $\ref{thm:main}$ are
 independent of the choice of coordinates of the source 
 and target.
 To be precise,
 the rank of\/ $(df)_p$,
 the non-degeneracy of\/ $p$,
 and the sign of\/ $\det\Hess\lambda(p)$,
 are independent of the choice of both coordinates on
 the source and target.
 Suppose further that $p$ is non-degenerate,
 and let $\lambda(u_1,u_2)$ and\/ $\tilde\lambda(v_1,v_2)$ are
 area density functions of\/ $f$,
 and\/ $\eta$ and\/ $\tilde\eta$ are
 null vector fields of\/ $f$,
 then the following hold:
 \begin{itemize}
  \item $\eta\lambda(p)=0$ if and only if\/
        $\tilde\eta\tilde\lambda(p)=0$.
  \item If\/ $\eta\lambda(p)=\tilde\eta\tilde\lambda(p)=0$, 
        then\/ $\eta\eta\lambda(p)\ne0$ if and only if\/
        $\tilde\eta\tilde\eta\tilde\lambda(p)\ne0$.
  \item If\/ $\eta\lambda(p)=\eta\eta\lambda(p)=
        \tilde\eta\tilde\lambda(p)=
        \tilde\eta\tilde\eta\tilde\lambda(p)=0$,
        then\/ $\eta\eta\eta\lambda(p)\ne0$ 
        if and only if\/
        $\tilde\eta\tilde\eta\tilde\eta\tilde\lambda(p)\ne0$.
 \end{itemize}
\end{lemma}
\begin{proof}
  Needless to say, $\rank(df)_p$
 is independent of the choice of the coordinate
  systems.
 If we change the coordinates,
 then the function $\lambda$ is multiplied by
 a non-zero function.
 Since the vanishing of
 $d\lambda$ and the sign of $\det\Hess\lambda$
 do not change under this multiplication,
 the first part of the lemma is proved.
 We now prove the second part.
 We can write $\tilde\eta=a_1\xi+a_2\eta$, where $a_1,a_2$
 are functions near $p$, satisfies $a_1=0$ on $S(f)$,
 and $\xi$ is a vector field transverse to $\eta$ at $p$,
 and assume that
 $\tilde\lambda$ is a multiplication of $\lambda$ by
 a non-zero function.
 Under this setting, since $\{\lambda=0\}=\{a_1=0\}$
 holds, one can prove that
 the non-degeneracy yields the 
 desired equivalences.
\end{proof}
Now we prove Theorem \ref{thm:main};
the method of proof is due to Rieger \cite{rie}.

\begin{proof}[Proof of $(1)$ and $(2)$.]
 Since $p$ is of corank one,
 $f$ can be represented
 as
 $$f(u,v)=(u,vf_2(u,v)),\quad p=(0,0)$$
 by Lemma \ref{lem:indeplamb}.
 Since $\lambda(p)=0$ and $d\lambda(p)=0$, we have $f_2=(f_2)_u=(f_2)_v=0$
 at $p$,
 where $(f_2)_u=\partial f_2/\partial u$ and
 $(f_2)_v=\partial
 f_2/\partial v$.
 Therefore, $f$ can be
 written as
 $$
 \big(u,v(au^2+2buv+cv^2)+\operatorname{(higher\ order\ term)}\big),
 \ a,b,c\in\R.
  $$
   Here, the
 ``higher order term'' consists of the terms
 whose degrees are greater than $3$.
 Since $\det\Hess\lambda(p)\ne0$,
 it holds that
 $a$, $b$ or $c$ does not vanish at $p$.
 Moreover, since $\eta=(0,1)$ and $\eta\eta\lambda(p)\ne0$,
 it holds that $c\ne0$.
 Now, by the coordinate change
 $$
 U=u,\quad V=v+\frac{2b}{3c}u,
 $$
 $f$ can be written as
 $$
 \big(u,v(\alpha u^2+\beta v^2)
 +
 \gamma u^3+\operatorname{(higher\ order\ term)}\big),\ 
 \alpha,\beta,\gamma\in\R.
 $$
 Here, the
 ``higher order term'' consists of the terms
 whose degrees are greater than $3$.
 We remark that the sign of 
 $\alpha\beta$ coincides with the sign of $\Hess\lambda(p)$.
 Hence, by some scaling change and a coordinate change on the target,
 $f$ can be written as
  \begin{equation}
   \label{eq:liptotyu}
  \big(u,v(u^2\pm v^2)+\operatorname{(higher\ order\ term)}\big).
  \end{equation}
  Since the map germ $(u,v(u^2\pm v^2))$
 is three-determined,
  the map germ (\ref{eq:liptotyu}) is ${\mathcal A}$-equivalent to the 
  lips$(+)$ or the beaks$(-)$.
\end{proof}
\begin{proof}[Proof of $(3)$.]
  Since $f$ is of corank one, $f$ can be written
  as $f(u,v)=(u,vh(u,v))$.
  Then the null vector field is $(0,1)$.
  Write
  $$
  \begin{array}{rcl}
  vh(u,v)&=&
  a_{11}uv+a_{02}v^2+a_{21}u^2v+a_{12}uv^2+a_{03}v^3
  +a_{31}u^3v\\
  &&\hspace{5mm}+a_{22}u^2v^2+a_{13}uv^3+a_{04}v^4
  +\operatorname{(higher\ order\ term)}.
  \end{array}
  $$
  Here, the ``higher order term'' consists of
  the terms whose degrees are greater than $4$.
  The non-degeneracy of $f$ yields that $a_{11}\ne0$.
  If $a_{02}\ne0$, by Fact \ref{fact:whit},
  $f$ is ${\mathcal A}$-equivalent to the fold.
  Moreover, if $a_{02}=0$ and $a_{03}\ne0$
  then by Fact \ref{fact:whit}, 
  $f$ is ${\mathcal A}$-equivalent to the cusp.
  Hence we can assume $a_{02}=a_{03}=0$.
 Since $\eta\eta\eta\lambda(p)\ne0$, we have $a_{04}\ne0$.
 By the coordinate change
 $$
 \begin{array}{l}
 \tilde u=u,\\
 \tilde v=a_{11}v+a_{21}uv+a_{12}v^2+a_{31}u^2v+a_{22}uv^2+a_{13}v^3,
 \end{array}
  $$
  $f$ is written as
  $$
  f(\tilde u,\tilde v)
  =
  \big(\tilde u,\tilde u\tilde v+\tilde v^4
  +
  \operatorname{(higher\ order\ term)}\big).
  $$
  Since $(\tilde u,\tilde u\tilde v+\tilde v^4)$ is four-determined,
  it is ${\mathcal A}$-equivalent to $(u,uv+v^4)$.
\end{proof}
\section{Singularities of characteristic surfaces of a single
conservation law}
\label{sec:cons}

In this section, we consider the
following Cauchy problem of a single conservation law:

\begin{equation}\tag{C}
 \label{singcons}
 \left\{
  \begin{array}{l}
   \dfrac{\partial y}{\partial t}(t,\x)
    +
    \sum_{i=1,2}\dfrac{d f_i}{d
    y}\big(y(t,\x)\big)\dfrac{\partial y}{\partial x_i}(t,\x)=0,\\
   y(0,\x)=\phi(\x),\quad \x=(x_1,x_2),
  \end{array}
 \right.
\end{equation}
where, $f_1,f_2$ and $\phi$ are functions.
We consider the characteristic surfaces of \eqref{singcons}
following the framework of \cite{ikcons}.

Let $\pi:PT^*(\R\times\R^2\times\R)\to\R\times\R^2\times\R$
be the projective cotangent bundle.
Identify $PT^*(\R\times\R^2\times\R)
=(\R\times\R^2\times\R)\times P(\R\times\R^2\times\R)$
and denote the local coordinates of this space by
$(t,\x,y,[\tau:\vect{\xi}:\eta])$.
We consider the canonical contact form,
$$
\alpha=[\tau dt+\xi_1dx_1+\xi_2dx_2+\eta dy].
$$
Then the equation \eqref{singcons} is written in the following
form:
$$
\begin{array}{l}
E(1,f_1',f_2',0)=
\Big\{\big(t,\x,y,[\tau:\vect{\xi}:\eta]\big)\in PT^*(\R\times\R^2\times\R)
\,\Big\vert\,\\
\hspace{75mm}\tau+\displaystyle\sum_{i=1,2}f_i'(y)\xi_i=0\Big\},
\end{array}
$$
where $f_i'=df_i/dy$.
If (\ref{singcons}) has a classical solution $y$,
then the non-zero normal vector $\nu=(y_t,y_{x_1},y_{x_2},-1)$
of smooth hypersurface
$(t,\x,y(t,\x))\subset\R\times\R^2\times\R$ exists,
where, $y_{x_1}=\partial y/\partial x_1$ for example.

Hence we have a Legendrian immersion 
$:\R\times\R^2\to PT^*(\R\times\R^2\times\R)$:
$$
\begin{array}{l}
\tilde{y}(t,\x):(t,\x)\mapsto
(t,\x,y(t,\x),[\nu])\\
\hspace{30mm}
\in
 E(1,f_{1}',f_{2}',0) \subset PT^*(\R\times\R^2\times\R).
\end{array}
$$
According to this, we define a {\it geometric solution\/} of
(\ref{singcons}) as a Legendrian immersion
$L:(U;u_1,u_2)\to E(1,f_{1}',f_{2}',0)\subset PT^*(\R\times\R^2\times\R)$
of a domain $U\subset\R^2$ such that $\pi\circ L$ is an embedding.
We apply the method of characteristic equation.
The characteristic equation associated with \eqref{singcons}
through $(0,\x_0)$ is
$$
\begin{array}{ll}
 \dfrac{dx_i}{dt}(t)=\dfrac{df_i}{dy}
 \Big(y\big(t,\x(t)\big)\Big),&\x(0)=\x_{0}\\
 \dfrac{dy}{dt}(t,\x(t))=0,&y(0,\x(0))=\phi(\x_0).
\end{array}
$$
The solution of the characteristic equation can be expressed
by
\begin{equation}
 \label{sol-chara}
 \begin{array}{l}
  x_i(\u,t)=u_i+t\dfrac{df_i}{dy}\big(\phi(\u)\big),\\[3mm]
  \hspace{20mm}
  y\big(0,\x(\u,0)\big)=y(0,\u)=\phi(\u),\quad  \u=(u_1,u_2)\in U.
  \end{array}
\end{equation}
If a map 
\begin{equation}
 \label{eq:charasurf}
  g_t:\u\mapsto \big(x_1(\u,t),x_2(\u,t)\big)
\end{equation}
is non-singular,
$y=\phi\big((g_t)^{-1}(x_1,x_2)\big)$
is the classical solution of (\ref{singcons})
(See \cite[Section 5]{ikcons}).
Remark that if $t=0$, $g_t$ is non-singular.
Thus, 
in order to investigate the singularity of \eqref{singcons},
we study the singularities of a family of maps $g_t$.
The discriminant function of $g_t(\u)$ is
$$
\det
\pmt{1+tc_{11}&tc_{12}\\
tc_{21}&1+tc_{22}},
\qquad
c_{ij}
=
\dfrac{d^2f_i}{dy^{2}}(\phi(\u))
\dfrac{\partial\phi}{\partial u_{j}}(\u).
$$
Needless to say, this matrix is never equal to the zero-matrix.
This implies that
$(t,\u)$ is a singular point of (\ref{sol-chara}),
if and only if $-t^{-1}$ is an eigen value of the matrix
$C=(c_{ij})_{i,j=1,2}$.
The eigen equation for an eigen value $\mu$ of $C$ can be
computed as
$$
\begin{array}{rcl}
\displaystyle
0&=&\det\left(C-\mu \pmt{1&0\\0&1}\right)\\[4mm]
&=&
\det
\pmt{
  \dfrac{d^2f_1}{dy^{2}}\big(\phi(\u)\big)
  \dfrac{\partial\phi}{\partial u_1}(\u)-\mu&
  \dfrac{d^2f_1}{dy^{2}}\big(\phi(\u)\big)
  \dfrac{\partial\phi}{\partial u_2}(\u)\\[4mm]
  \dfrac{d^2f_2}{dy^{2}}\big(\phi(\u)\big)
  \dfrac{\partial\phi}{\partial u_1}(\u)&
  \dfrac{d^2f_2}{dy^{2}}\big(\phi(\u)\big)
  \dfrac{\partial\phi}{\partial u_2}(\u)-\mu}
\\[11mm]
&=&
\mu\left(
\mu-\dfrac{d^2f_1}{dy^{2}}\big(\phi(\u)\big)
\dfrac{\partial\phi}{\partial u_1}(\u)-
\dfrac{d^2f_2}{dy^{2}}\big(\phi(\u)\big)
\dfrac{\partial\phi}{\partial u_2}(\u)\right)\\[4mm]
&=&
\mu(\mu-\trace C).
\end{array}
$$
Hence $(t,\u)$ is a singular point of (\ref{sol-chara}),
if and only if
$$
t=-1/\trace{C}.
$$
We call $C$ the {\it shape operator\/} of (\ref{singcons}).
Now we consider the first singular point of (\ref{eq:charasurf})
with respect to $t$ from the initial time $t=0$.

For a minimal value of
$t(\u)=-1/\trace{C}$,
if $\det\Hess{t(\u)}>0$ holds,
then by Theorem \ref{thm:main}, the singular point at $\u$ 
is ${\mathcal A}$-equivalent to the lips.
Izumiya and Kossioris \cite{ikcons}
have developed an unfolding theory and
classified the generic singularities
of multi-valued solutions in general dimensions.
According to it,
the first singular point of (\ref{eq:charasurf}) is
generically the lips,
where they did not give a
condition for the singular point to be
equivalent to the lips.
Using our criterion for the lips,
we detect the singular point and
write down an explicit condition for
the singular point to be equivalent to
the lips.
As a corollary of it, we give a simple proof that
the first singular point of (\ref{eq:charasurf}) is
generically the lips.

Since the single conservation law (\ref{singcons}) is
determined by functions $(f_1,f_2)$ 
and the initial value $\phi$,
we may regard that the space of 
single conservation laws is
the space
$$
\{(f_1,f_2,\phi)\}=
C^\infty(\R,\R)^2\times C^\infty(\R^2,\R)
$$
with the Whitney $C^\infty$-topology.
\begin{theorem}
 There exists a residual subset\/
 ${\mathcal O}\subset C^\infty(\R,\R)^2\times C^\infty(\R^2,\R)$
 such that for any\/ $(f_1,f_2,\phi)\in{\mathcal O}$,
 the map germ\/ $(\ref{eq:charasurf})$ defined by\/ $(f_1,f_2,\phi)$ at
 the first singular point with respect to\/ $t>0$
 is\/ ${\mathcal A}$-equivalent to the lips.
\end{theorem}
Here, a subset is {\it residual\/} if it is
a countable intersection of open and dense subsets.
\begin{proof}
 Since, for a function $w$, the behaviors of $dw$ and $\Hess w$ are the
 same as those of $w^{-1}$,
 we may calculate these quantities about $\trace C$.
 By a direct calculation, we have
 $$
 \begin{array}{l}
  \Xi_1(\u):=(1/t)_{u_1}=
   f_1^{(3)}\,
   (\phi_{1})^2 + 
   f_2^{(3)}\,
   \phi_{1}\,\phi_{2} + 
   f_1''\,\phi_{11}+   f_2''\,
   \phi_{12},\\
  \Xi_2(\u):=  (1/t)_{u_2}=
   f_1^{(3)}\,
   \phi_{1}\,\phi_{2}+f_2^{(3)}\,
   (\phi_{2})^2 + 
   f_1''\,\phi_{12}
   +   f_2''\,
   \phi_{22}
   \end{array}
 $$
 and
 $$
 \begin{array}{l}
  \Xi_3(\u):=\det\Hess{(1/t)}=\\
  f_1^{(4)}\Bigg[
    f_1^{(3)}\,(\phi_{1})^2\Big(
      (\phi_{1})^2\,\phi_{22}
      -2\,\phi_{1}\,\phi_{2}\,\phi_{12}
      +(\phi_{2})^2\,\phi_{11}
    \Big)\\
    \hspace{14mm}
    +f_1''\,\phi_{1}\Big((\phi_{1})^2\,\phi_{122}
    -2\,\phi_{1}\,\phi_{2}\,\phi_{112}
    +(\phi_{2})^2\,\phi_{111}\Big)\\\hspace{16mm}
    +f_2^{(3)}\,\phi_{1}\,\phi_{2}\,
    \Big((\phi_{1})^2\,\phi_{22}
       - 
      2\,\phi_{1}\,\phi_{2}\,
      \phi_{12} + 
      (\phi_{2})^2\,\phi_{11}
    \Big)\\\hspace{22mm}
    +f_2''\,\phi_{1}\,
      \Big((\phi_{1})^2\, \phi_{222}
       - 
      2\,\phi_{1}\,\phi_{2}\,
      \phi_{122} + 
      (\phi_{2})^2\,\phi_{112}
    \Big)
  \Bigg]
\\
  +f_2^{(4)}\Bigg[
    f_1^{(3)}\phi_{1}\phi_{2}\Big(
      (\phi_{1})^2\,\phi_{22}
      -2\,\phi_{1}\,\phi_{2}\,\phi_{12}
      +(\phi_{2})^2\,\phi_{11}
    \Big)\\
    \hspace{14mm}
    +f_1''\phi_{2}\Big(
      (\phi_{1})^2\,\phi_{122}
      -2\,\phi_{1}\,\phi_{2}\,\phi_{112}
      +(\phi_{2})^2\,\phi_{111}
    \Big)
    \\\hspace{16mm}
    +f_2^{(3)}\,(\phi_{2})^2\,\Big(
      (\phi_{1})^2\,\phi_{22}
      -2\,\phi_{1}\,\phi_{2}\,
      \phi_{12} + 
      (\phi_{2})^2\,\phi_{11}
    \Big)\\\hspace{22mm}
    +f_2''\,\phi_{2}\,\Big(
      (\phi_{1})^2\,\phi_{222}      
      -2\,\phi_{1}\,\phi_{2}\,
      \phi_{122} + 
      (\phi_{2})^2\,\phi_{112}
    \Big)
  \Bigg]
\\
  +(f_1^{(3)})^2\Big[
     \phi_{1}\phi_{11}(3\phi_{1}\phi_{22}
     +2\phi_{2}\phi_{12})-4(\phi_{1})^2(\phi_{12})^2
     -(\phi_{2})^2(\phi_{11})^2\Big]
\\[2mm]
  + 
  (f_2^{(3)})^2\Big[
    \phi_{2}\phi_{22}
    (2\phi_{1}\phi_{12} + 3\phi_{2}\phi_{11})
    -(\phi_{1})^2(\phi_{22})^2
    -4(\phi_{2})^2(\phi_{12})^2
  \Big]
\\[2mm]  +
  f_1^{(3)}f_2^{(3)}
    \Big[
      -2\,(\phi_{1})^2\,\phi_{12}\,\phi_{22}
      -4\,\phi_{1}\,\phi_{2}\,(\phi_{12})^2\\
      \hspace{35mm}
      +8\,\phi_{1}\,\phi_{2}\,\phi_{11}\,\phi_{22}
      -2\,(\phi_{2})^2\,\phi_{11}\,\phi_{12}
    \Big]
\\
  +f_1^{(3)}
  \Bigg[
    f_1''\Big(
      3\,\phi_{1}\,\phi_{11}\,\phi_{122}
      -4\,\phi_{1}\,\phi_{12}\,\phi_{112}
      \\\hspace{30mm}
      -2\,\phi_{2}\,\phi_{11}\,\phi_{112}
      +\phi_{1}\,\phi_{22}\,\phi_{111}
      +2\,\phi_{2}\,\phi_{12}\,\phi_{111}
    \Big)
    \\\hspace{12mm}
    +f_2''\Big(
      -4\phi_{1}\,\phi_{12}\,\phi_{122}
      +3\,\phi_{1}\,\phi_{11}\,\phi_{222}
      \\\hspace{30mm}
      -2\,\phi_{2}\,\phi_{11}\,\phi_{122}
      +\phi_{1}\,\phi_{22}\,\phi_{112}
      +2\,\phi_{2}\,\phi_{12}\,\phi_{112}
    \Big)
  \Bigg]
      \end{array}
  $$
  $$
  \begin{array}{l}
  +f_2^{(3)}\,
  \Bigg[
    f_1''\Big(
      2\phi_{1}\,\phi_{12}\,\phi_{122}
      +\phi_{2}\,\phi_{11}\,\phi_{122}
      \\\hspace{30mm}
      -2\,\phi_{1}\,\phi_{22}\,\phi_{112}
      -4\,\phi_{2}\,\phi_{12}\,\phi_{112}
      +3\,\phi_{2}\,\phi_{22}\,\phi_{111}
    \Big)\\\hspace{12mm}
    +f_2''
    \Big(
      2\phi_{1}\,\phi_{12}\,\phi_{222}
      -2\,\phi_{1}\,\phi_{22}\,\phi_{122} 
      \\\hspace{30mm}
      -4\,\phi_{2}\,\phi_{12}\,\phi_{122}
      +\phi_{2}\,\phi_{11}\,\phi_{222}
      +3\,\phi_{2}\,\phi_{22}\,\phi_{112}
    \Big)
  \Bigg]\\
    +
    (f_1'')^2\,
    \Big(\phi_{111}\,\phi_{122}-(\phi_{112})^2\Big)
    +(f_2'')^2\,
    \Big(\phi_{112}\,\phi_{222} -(\phi_{122})^2\Big)\\
      \hspace{60mm}+ 
   f_1''f_2''\,
   \Big( 
      \phi_{111}\,\phi_{222}-\phi_{112}\,\phi_{122}\Big),
 \end{array}
 $$
 where for the sake of simplicity, we set
 $$
\begin{array}{l}
 f_{\ell}':=\dfrac{d f_{\ell}}{dy},\ 
  f_{\ell}'':=\dfrac{d^2 f_{\ell}}{dy^2},\
  f_{\ell}^{(m)}:=\dfrac{d^m f_{\ell}}{dy^m},\
  (\ell=1,2,\ m=3,4)
\\
  \phi_{i}:=\dfrac{\partial \phi}{\partial u_i},\ 
 \phi_{ij}:=\dfrac{\partial^2 \phi}{\partial u_iu_j},\
\text{and}\ 
 \phi_{ijk}:=
 \dfrac{\partial^3 \phi}{\partial u_iu_ju_k},\ 
 (i,j,k=1,2).
\end{array} 
$$
 Next we consider a map
 $$
 \begin{array}{l}
  j^4(f_1,f_2,\phi):(y,\u)\mapsto
   \big(j^4f_1(y),j^4f_2(y),j^4\phi(\u)
   \big)\\
   \hspace{50mm}\in
 J^4(\R,\R)^2\times J^4(\R^2,\R)
 \end{array}
 $$
 and four subsets of jet spaces $J^4(\R,\R)^2\times J^4(\R^2,\R)$
 as follows:
 $$
 \begin{array}{l}
  \hat{\Xi}_0:=\big\{j^4(f_1,f_2,\phi)(y,\u)
   \,\vert\,y-\phi(\u)=0\big\}\\
  \hat{\Xi}_1:=\big\{j^4(f_1,f_2,\phi)(y,\u)
   \,\vert\,\Xi_1(\u)=0\big\}\\
  \hat{\Xi}_2:=\big\{j^4(f_1,f_2,\phi)(y,\u)
   \,\vert\,\Xi_2(\u)=0\big\}\\
  \hat{\Xi}_3:=\big\{j^4(f_1,f_2,\phi)(y,\u)
   \,\vert\,\Xi_3(\u)=0\big\}.
 \end{array}
 $$
 Since the coordinate system of $J^4(\R,\R)^2\times J^4(\R^2,\R)$
 is defined by each coordinate of source and
 value of derivatives of functions,
 $\hat\Xi_0$, $\hat\Xi_1$, $\hat\Xi_2$ and $\hat\Xi_3$
  are algebraic subsets with respect to
 the coordinates of $J^4(\R,\R)^2\times J^4(\R^2,\R)$.
 Comparing the coefficients of
 $\phi_{11}$ and $\phi_{22}$ in $\Xi_1$ and $\Xi_2$,
 we see that $\Xi_1$ and $\Xi_2$ do not have a common factor.
 Moreover,
 $f_1''\phi_{1}(\phi_{2})^2$ is the coefficient of
 $\phi_{111}f_1^{(4)}$ of $\Xi_3$, but
 this does not appear in either $\Xi_1$ or $\Xi_2$.
 Hence $S:=\cap_{i=0}^3\hat{\Xi}_i$ is
 a closed algebraic subset with codimension $4$
 in $J^4(\R,\R)^2\times J^4(\R^2,\R)$.
 So this set has a standard stratification.
 Applying the Thom jet transversality theorem to
 $j^4(f_1,f_2,\phi)$ and $S$,
 there exists a residual subset
 ${\mathcal O}\subset
 C^\infty(\R,\R)^2\times C^\infty(\R^2, \R)$
 such that for any $(f_1,f_2,\phi)\in {\mathcal O}$, the map
 $j^4(f_1,f_2,\phi)$ is transverse to $S$.
 Since the codimension of $S$ is 4,
 transversal condition means having no intersection point.
 If $g_{t_0}$ at $\u_0$ is the beaks, there is a singularity
 of $g_{t-\ep}$ near $\u_0$ for a sufficiently small number $\ep$,
 the beaks never appear at the minimal value of $t$ which $g_t$ is
 singular.
 Thus $(f_1,f_2,\phi)\in {\mathcal O}$ satisfies the desired condition. 
\end{proof}
\smallskip

The author would like to thank Professors
Goro Akagi, Shyuichi Izumiya and Farid Tari
for fruitful discussions and helpful advices.

\medskip

\begin{flushright}
\begin{tabular}{l}
Department of Mathematics,\\
Faculty of Education,\\
Gifu University,\\
Yanagido 1-1, Gifu, 501-1193, Japan.\\
e-mail: {\tt ksajiO\!\!\!agifu-u.ac.jp}
 \end{tabular}\end{flushright}
\end{document}